\documentclass[letterpaper,11pt,oneside]{amsart}
\usepackage{amsmath}
\usepackage{amssymb}
\usepackage{latexsym}
\usepackage{amsfonts}
\usepackage[latin1]{inputenc}
\usepackage[english]{babel}
\usepackage{multicol}
\usepackage{oldgerm}
\usepackage{enumerate}
\newtheorem{theorem}{Theorem}
\newtheorem{corollary}{Corollary}
\newtheorem{lemma}{Lemma}

\begin{document}

\title{Coincidences in generalized Lucas sequences}

\author{Eric F. Bravo}
\address{Departamento de Matem\'aticas\\ Universidad del Cauca\\ Calle 5 No 4--70\\Popay\'an, Colombia.}
\email{fbravo@unicauca.edu.co}

\author{Jhon J. Bravo}
\address{Departamento de Matem\'aticas\\ Universidad del Cauca\\ Calle 5 No 4--70\\Popay\'an, Colombia.}
\email{jbravo@unicauca.edu.co}

\author{Florian Luca}
\address{Mathematical Institute, UNAM Juriquilla\\  76230 Santiago de Quer\'etaro, M\'exico\\  and\\ 
School of Mathematics, University of the Witwatersrand\\  P. O. Box Wits 2050, South Africa}
\email{fluca@matmor.unam.mx}


\date{\today}

\begin{abstract}
For an integer $k\geq 2$,  let $(L_{n}^{(k)})_{n}$ be the $k-$generalized Lucas sequence which starts with $0,\ldots,0,2,1$ ($k$ terms) and each term afterwards is the sum of the $k$ preceding  terms. In this paper, we find all the integers that appear in different generalized Lucas sequences; i.e., we study the Diophantine equation $L_n^{(k)}=L_m^{(\ell)}$ in nonnegative integers $n,k,m,\ell$ with $k, \ell\geq 2$. The proof of our main theorem uses lower bounds for linear forms in logarithms of algebraic numbers and a version of the Baker--Davenport reduction method.  This paper is a continuation of the earlier work \cite{BL3}.  

\medskip



\end{abstract}

\maketitle


\section{Introduction}

\noindent Let $k\geq 2$ be an integer. We consider the linear recurrence sequence of order $k$ denoted $G^{(k)}:=(G_n^{(k)})_{n\geq 2-k}$ defined as 
\begin{equation*}\label{recurrenciaG}
G_{n}^{(k)}=G_{n-1}^{(k)}+G_{n-2}^{(k)}+\cdots+G_{n-k}^{(k)}\quad \text{for all} \quad n\ge 2,
\end{equation*}
with the initial conditions $G^{(k)}_{-(k-2)}=G^{(k)}_{-(k-3)}=\cdots=G^{(k)}_{-1}=0$, $G^{(k)}_0=a$ and $G^{(k)}_1=b$.  

Observe that if $a=0$ and $b=1$, then $G^{(k)}$ is nothing more than the $k-$Fibonacci sequence $F^{(k)}:=(F_n^{(k)})_{n\geq 2-k}$. In this case, if we choose $k=2$, we obtain the classical Fibonacci sequence $(F_n)_{n\geq 0}$. On the other hand, if $a=2$ and $b=1$, then $G^{(k)}$ is known as the $k-$Lucas sequence $L^{(k)}:=(L_n^{(k)})_{n\geq 2-k}$. In the special case of $k=2$, we obtain  the usual Lucas companion of the Fibonacci sequence 
\[
L_0=2, \quad L_1=1 \quad\text{and}\quad L_{n}=L_{n-1}+L_{n-2} \quad \text{for} \quad n\geq 2.
\]
\[
(L_{n})_{n\ge 0}=\{2, 1, 3, 4, 7, 11, 18, 29, 47, 76, 123, 199, 322, 521, 843, 1364,\ldots\}.
\]
For example, if $k=3$, then the $3-$Lucas sequence is
\[
(L_{n}^{(3)})_{n\ge -1}=\{0, 2, 1, 3, 6, 10, 19, 35, 64, 118, 217, 399, 734, 1350, 2483, 4567,\ldots\}.
\]
If $k=4$, we get the $4-$Lucas sequence
\[
(L_{n}^{(4)})_{n\ge -2}=\{0, 0, 2, 1, 3, 6, 12, 22, 43, 83, 160, 308, 594, 1145, 2207, 4254, 8200,\ldots\}.
\]
As can be seen in \cite[Lemma 2]{BLIndia}, these generalized Lucas sequences have the remarkable property that  the first few terms are given by 
\[
L_{n}^{(k)}=3\cdot 2^{n-2} \quad \text{for all} \quad 2\leq n \leq k.
\]
The above sequences are among the several generalizations of the Fibonacci numbers which have been studied in the literature.  Other 
generalizations are also known (see, for example, \cite{Brent,Kilic,Muskat}).

Several authors have worked on problems involving generalized Fibonacci sequences. For instance, F. Luca \cite{FL} and D. Marques \cite{DM} proved that 55 and 44 are the largest repdigits in the sequences $F^{(2)}$ and $F^{(3)}$, respectively. Moreover, D. Marques conjectured that there are no repdigits with at least two digits in $F^{(k)}$ for any $k > 3$. This conjecture was confirmed in \cite{BL1}. In addition, the Diophantine equation $F_n^{(k)}=2^m$ was studied in \cite{BL2}. 

In 2005, T. D. Noe and J. V. Post \cite{TonyNoe} proposed a conjecture about coincidences of terms of generalized Fibonacci sequences. In their work, they gave a heuristic argument to show that if $k\neq \ell$, then the cardinality of the intersection $F^{(k)}\cap F^{(\ell)}$ must be small. Further, they used computational methods which led them to confirm the conjecture for all terms whose magnitude is less than 22000. This conjecture has been recently proved to hold independently by Bravo--Luca \cite{BL3} and D. Marques \cite{DM1}. 

In this paper, we investigate the problem of determining the intersection of two generalized Lucas sequences. To begin with, it is important to mention that Mignotte (see \cite{Mignotte}) proved (under some technical conditions) that only a finite number of coincidences between two fixed linear recurrence sequences can occur. In this context, one could of course ask how large is the cardinality of the finite set $L^{(k)}\cap L^{(\ell)}$ for $k>\ell\ge 2$.
From the above initial values, we see that there are some numbers that appear in different generalized Lucas sequences. For instance, the zeros that appear at the beginning, but these numbers are not interesting for us. Throughout this paper we only consider nonzero terms of these sequences.

Here, we determine all the solutions of the Diophantine equation
\begin{equation} \label{eqn}
L_n^{(k)}=L_m^{(\ell)}, 
\end{equation}
in nonnegative integers $n,k,m,\ell$ with $k>\ell\geq 2$.

First of all, note that if $k>\ell$, then $L_t^{(k)}=L_t^{(\ell)}$ for all $0\leq t \leq \ell$, i.e., the quadruple
\begin{equation}\label{trivial}
(n,k,m,\ell)=(t,k,t,\ell),
\end{equation}
is a solution of equation \eqref{eqn} for all $0\leq t \leq \ell$. Solutions given by \eqref{trivial}  will be called \emph{trivial solutions}.  

We prove the following theorem.

\begin{theorem} \label{teo1}
The Diophantine equation \eqref{eqn} has only trivial solutions. 
\end{theorem}
 As immediate consequences of Theorem \ref{teo1} we have the following corollaries.
\begin{corollary}
Let $k,\ell$ be integers with $k>\ell\geq 2$. Then 
\[
|L^{(k)}\cap L^{(\ell)}|=\ell+1.
\]
\end{corollary}
\begin{corollary}
If $(n,k,a)$ is a solution of the Diophantine equation $L_{n}^{(k)}=3\cdot 2^{a}$  in nonnegative integers $n,k,a$ with $k\geq 2$, then $0\leq n \leq k$ and $a=n-2$.
\end{corollary}

In this paper, we follow the approach and the presentation described in \cite{BL3}.


\section{Preliminary results}

\noindent Before proceeding further, we recall some facts and properties of the $k$-generalized Lucas sequences which will be used later. First, it is known that the characteristic polynomial of the sequence $G^{(k)}$, namely 
\[
\Psi_k(x)=x^k-x^{k-1}-\cdots-x-1,
\]
is irreducible over $\mathbb{Q}[x]$ and has just one root outside the unit circle; the other roots are strictly inside the unit circle (see, for example, \cite{Mi}, \cite{Mil} and \cite{DAW}). Throughout this paper, $\alpha:=\alpha(k)$ denotes that single root, which is located between $2(1-2^{-k})$ and 2 (see \cite{DAW}).  We shall use $\alpha_1,\ldots,\alpha_k$ for all the roots of  $\Psi_k(x)$ with the convention that $\alpha_1:=\alpha$. Similarly, we use $\beta_1,\ldots,\beta_{\ell}$ for the roots of $\Psi_\ell(x)$, with the convention that $\beta_1:=\beta$ is the real root of $\Psi_\ell(x)$ exceeding $1$.

We now consider for an integer $s\geq 2$, the function 
\begin{equation}\label{fun-f}
f_s(x)=\frac{x-1}{2+(s+1)(x-2)} \quad \text{for} \quad x>2(1-2^{-s}). 
\end{equation}
With this notation, the following ``Binet--like" formula for $F^{(k)}$ appears in  Dresden \cite{GPD}:
\[
F_n^{(k)}=\sum_{i=1}^{k}f_k(\alpha_i)\alpha_i^{n-1}.
\]
It was also proved in \cite{GPD} that the approximation 
\[
|F_n^{(k)}-f_k(\alpha)\alpha^{n-1}|<\frac{1}{2}  \quad \text{holds~for all}\quad n\geq 2-k.
\]
Further, in \cite{BL1}, it is proved that
\[
\alpha^{n-2} \leq F_n^{(k)}\leq \alpha^{n-1} \quad \text{for all}\quad n\geq 1\quad \text{and}\quad k\ge 2.
\]
Analogous results to the previous facts have recently been established by Bravo and Luca \cite{BLIndia} for the sequence $L^{(k)}$.

\begin{lemma}[Properties of $L^{(k)}$] \label{properties-L}
Let $k\geq 2$ be an integer. Then
\begin{itemize}
\item[$(a)$]  $\alpha^{n-1} \leq L_n^{(k)}\leq 2\alpha^{n}$ for all $n\geq 1$.
\item[$(b)$] $L^{(k)}$ satisfies the following ``Binet--like" formula
      \[
         L_n^{(k)}=\sum_{i=1}^{k}(2\alpha_i-1)f_k(\alpha_i)\alpha_i^{n-1},
      \]     
     where $\alpha=\alpha_1,\ldots,\alpha_k$ are the roots of $\Psi_k(x)$.
\item[$(c)$] $|L_n^{(k)}-(2\alpha-1)f_k(\alpha)\alpha^{n-1}|<3/2$  holds for all $n\geq 2-k$.
\end{itemize}
\end{lemma}

Now assume that we have a nontrivial solution $(n,k,m,\ell)$ of equation \eqref{eqn} with the previous conventions that $\alpha=\alpha(k)$ and $\beta=\alpha(\ell)$. By Lemma \ref{properties-L} $(a)$, we have
\[
\beta^{m-1} \leq L_m^{(\ell)}=L_n^{(k)}\leq 2\alpha^{n}<2^{n+1},
\]
so, we get
\begin{equation}\label{util2}
m <\frac{3n+5}{2}, \quad\text{or, equivalently}\quad \frac{2m-5}{3}<n,
\end{equation}
where we have used the fact that the inequality $1/\log\beta<2.1$ holds for all $\ell\geq 2$. We record this estimate for future referencing.

In order to prove Theorem \ref{teo1}, we need to use several times a Baker--type lower bound for a nonzero linear form in logarithms of algebraic numbers and such a bound, which plays an important role in this paper, was given by Matveev \cite{Matveev}. We begin by recalling some basic notions from algebraic number theory.

Let $\eta$ be an algebraic number of degree $d$ with minimal polynomial over the integers
\[
a_0x^d+a_1x^{d-1}+\cdots+a_d=a_0\prod_{i=1}^{d}(X-\eta^{(i)}),
\]
where the $a_i$'s are relatively prime integers with $a_0>0$ and the $\eta^{(i)}$'s are conjugates of $\eta$. Then
\[
h(\eta)=\frac{1}{d}\left(\log a_0+\sum_{i=1}^{d}\log\left(\max\{|\eta^{(i)}|,1\}\right)\right)
\]
is called the \emph{logarithmic height} of $\eta$. In particular, if $\eta=p/q$ is a rational number with $\gcd(p,q)=1$ and $q>0$, then $h(\eta)=\log \max \{|p|,q\}$. 

The following properties of the logarithmic heigh, which will be used in the next sections without special reference, are also known:
\begin{itemize}
\item $h(\eta\pm\gamma)\leq h(\eta) + h(\gamma)+\log 2$.
\item $h(\eta\gamma^{\pm 1})\leq h(\eta)+h(\gamma)$.
\item $h(\eta^{s})=|s|h(\eta)$.
\end{itemize}

With the previous notation, Matveev (see \cite{Matveev} or Theorem 9.4 in \cite{Bug}) proved the following deep theorem.

\begin{theorem}[Matveev's theorem]\label{teoMatveev}
Assume that $\gamma_1, \ldots, \gamma_t$ are positive real algebraic numbers in a real algebraic number field $\mathbb{K}$ of degree $D$, $b_1,\ldots,b_t$ are rational integers, and 
\[
\Lambda:=\gamma_1^{b_1}\cdots\gamma_t^{b_t}-1,
\]
is not zero. Then
\[
|\Lambda|>\exp\left(-1.4\times 30^{t+3}\times t^{4.5}\times D^2(1+\log D)(1+\log B)A_1\cdots A_t\right),
\]
where
\[
B\geq \max\{|b_1|,\ldots,|b_t|\},
\]
and
\[
A_i\geq \max\{Dh(\gamma_i),|\log \gamma_i|, 0.16\}, \quad \text{for all} \quad  i=1,\ldots,t.
\] 
\end{theorem}

We will also use the following estimates from \cite{BLIndia}. A key point of that work consists of exploiting the fact that when $k$ is large, the dominant root of $L^{(k)}$ is exponentially close to 2, so one can write the dominant term of the Binet formula for $L^{(k)}$ as 3 times a power of 2 plus an error which is well under control. Let us state this result as a lemma since we have some use for it later.

\begin{lemma} \label{estimaciondeltaeta}
For $k\geq 2$, let $\alpha$ be the dominant root of the characteristic polynomial $\Psi_k(x)$ of the $k-$Lucas sequence, and consider the function $f_k(x)$ defined in \eqref{fun-f}. Then
\[
h((2\alpha-1)f_k(\alpha))<\log 3+3\log k,
\]
where $h(\cdot)$ represents the logarithmic height function. Moreover, if $r>1$ is an integer satisfying $r-1<2^{k/2}$, then  
\[
(2\alpha-1)f_k(\alpha)\alpha^{r-1}=3\cdot 2^{r-2}+3\cdot 2^{r-1}\eta+\frac{\delta}{2}+\eta\delta,
\]
where $\delta$ and $\eta$ are real numbers such that
\[
|\delta|<\frac{2^{r+2}}{2^{k/2}}\quad  \text{and} \quad |\eta|<\frac{2k}{2^{k}}.
\]
\end{lemma}

In 1998, Dujella and Peth\H o in \cite[Lemma 5$(a)$]{DP} gave a version of the reduction method based on the Baker--Davenport lemma \cite{Baker-Davenport}.  We next present the following lemma from \cite{BL1}, which is an immediate variation of the result due to Dujella and Peth\H o from \cite{DP}, and will be one of the key tools used in this paper to reduce the upper bounds on the variables of the Diophantine equation \eqref{eqn}.

\begin{lemma} \label{reduce}
Let $M$ be a positive integer, let $p/q$ be a convergent of the continued fraction of the irrational $\gamma$ such that $q>6M$, and let $A,B,\mu$ be some real numbers with $A>0$ and $B>1$. Let $\epsilon:=||\mu q||-M||\gamma q||$, where $||\cdot||$ denotes the distance from the nearest integer. If $\epsilon >0$, then there is no solution to the inequality
\[
0<u\gamma-v+\mu<AB^{-w},
\]
in positive integers $u,v$ and $w$ with
\[
u\leq M \quad\text{and}\quad w\geq \frac{\log(Aq/\epsilon)}{\log B}.
\]
\end{lemma}


\section{An inequality for $n$ and $m$ in terms of $k$}

\noindent Since $k>\ell$ and the solution to equation \eqref{eqn} is nontrivial, we get easily that $m>n\geq 6$. Thus, in the remainder of the article, we can suppose that $\ell \leq m-1$, for otherwise there is nothing to prove.   

We now argue as in \cite{BL3}. Indeed, by using \eqref{eqn} and Lemma \ref{properties-L} $(c)$, we get that 
\begin{equation}\label{des-al-beta}
\begin{split}
|(2\alpha-1)f_k(\alpha)&\alpha^{n-1}-(2\beta-1)f_{\ell}(\beta)\beta^{m-1}|\\
& =|((2\alpha-1)f_k(\alpha)\alpha^{n-1}-L_n^{(k)})+(L_m^{(\ell)}-(2\beta-1)f_{\ell}(\beta)\beta^{m-1})|< 3.
\end{split}
\end{equation}
Dividing both sides of the above inequality by $(2\beta-1)f_{\ell}(\beta)\beta^{m-1}$, which is positive, we obtain
\begin{equation} \label{deslambda}
\left|\alpha^{n-1}\cdot\beta^{-(m-1)}\cdot (2\alpha-1)f_k(\alpha)((2\beta-1)f_{\ell}(\beta))^{-1}-1\right|<\frac{6}{\beta^{m-1}},
\end{equation}
where we used the fact that $1/f_{\ell}(\beta)<4$, which is easily seen taking into account that $2+(\ell+1)(\beta-2)<2$ and $1/(\beta-1)<2$. 

In a first application of Matveev's theorem, we take $t:=3$ and
\[
\gamma_1:=\alpha, \quad \gamma_2:=\beta, \quad \gamma_3:=(2\alpha-1)f_k(\alpha)((2\beta-1)f_{\ell}(\beta))^{-1}.
\]
We also take $b_1:=n-1$, $b_2:=-(m-1)$ and $b_3:=1$. Hence, 
\[
\Lambda:=\gamma_1^{b_1}\cdot\gamma_2^{b_2}\cdot\gamma_3^{b_3}-1.
\]
The algebraic number field containing $\gamma_1,\gamma_2,\gamma_3$ is $\mathbb{K}:=\mathbb{Q}(\alpha,\beta)$. Then $D=[\mathbb{K}:\mathbb{Q}]\leq k\ell$. The proof that $\Lambda\neq 0$ is similar to that given in \cite[p. 2126]{BL3}. We include it here for the sake of completeness.

Arguing by contradiction let us assume that $\Lambda=0$. Then
\begin{equation}\label{lambdanozero}
\frac{(2\alpha-1)(\alpha-1)}{2+(k+1)(\alpha-2)}\alpha^{n-1}=\frac{(2\beta-1)(\beta-1)}{2+(\ell+1)(\beta-2)}\beta^{m-1}.
\end{equation}
Let ${\mathbb L}={\mathbb Q}(\alpha_1,\ldots,\alpha_k,\beta_1,\ldots,\beta_\ell)$ be the normal closure of ${\mathbb K}$ and let further $\sigma_1,\ldots,\sigma_k$ be elements of ${\text{\rm Gal}}({\mathbb L}/{\mathbb Q})$ such that $\sigma_i(\alpha)=\alpha_i$. Since $k>\ell$, there exist $i\ne j$ in $\{1,2,\ldots,k\}$ such that $\sigma_i(\beta)=\sigma_j(\beta)$. Applying $\sigma_j^{-1}\sigma_i$ to the relation \eqref{lambdanozero} and then taking absolute values, we get that
\begin{equation} \label{lambdanozero1}
\left|\frac{(2\alpha_s-1)(\alpha_s-1)}{2+(k+1)(\alpha_s-2)} \alpha_s^{n-1}\right|=\frac{(2\beta-1)(\beta-1)}{2+(\ell+1)(\beta-2)}\beta^{m-1}
\end{equation}
where $s\neq 1$ is such that $\sigma_j^{-1}(\alpha_i)=\alpha_s$. But the above relation \eqref{lambdanozero1} is not possible since its left--hand side is smaller than 3, because $|\alpha_s|<1$ and
\[
|2+(k+1)(\alpha_s-2)|\geq (k+1)|\alpha_s-2|-2> k-1\geq 2,
\]
while  its right--hand side exceeds $L_m^{(\ell)}-3/2 > 3$ since $m\geq 7$.  Thus, $\Lambda\neq 0$.

Since $h(\gamma_1)=(\log\alpha)/k<(\log 2)/k=(0.693147\ldots)/k$ and $D\leq k\ell$, it follows that we can take 
$A_1:=0.7k>0.7\ell> Dh(\gamma_1)$. Similarly, we can take $A_2:=0.7k$.

We now observe that, by Lemma \ref{estimaciondeltaeta}, we have that
\[
h(\gamma_3)\leq \log 9+3\log k+3\log \ell<\log 9+6\log k\leq 8\log k
\]
for all $k\geq 3$. So, we can take $A_3:=8k^2\log k$. By recalling that $n<m$, we can take $B:=m-1$. Applying Theorem \ref{teoMatveev} to get a lower bound for $|\Lambda|$ and comparing it with inequality \eqref{deslambda}, we get
\[
\exp\left(-C_1(k)\times(1+\log (m-1))\,(0.7k)\,(0.7k)\,(8k^2\log k)\right)<\frac{6}{\beta^{m-1}},
\]
where $C_1(k):=1.4\times 30^{6}\times 3^{4.5}\times D^2 \times(1+\log D)<1.5 \times 10^{11}\,k^4\,(1+2\log k)$. Taking logarithms on both sides and performing the respective calculations, we get that
\begin{equation}\label{deslog}
\frac{m-1}{\log(m-1)}<7.41\times 10^{12}\,k^8\,\log^2k,
\end{equation}
giving
\[
m-1<5.34\times 10^14 \,k^8\,\log^2k.
\]
In the above we used the fact that inequality $x\log x<A$ implies $x<2A\log A$ whenever $A\geq 3$ (see \cite[p. 74]{BL2}). Let us record this result for future use. 

\begin{lemma}\label{cota_nl}
If $(n,k,m,\ell)$ is a nontrivial solution in positive integers of equation \eqref{eqn} with $k>\ell\ge 2$, then $\ell\leq m-1$ and the inequalities  
\[
6\leq n<m<5.4\times 10^{14}\,k^8\,\log^3 k
\]
hold.
\end{lemma}


\section{The case of small $k$}

\noindent We next treat the cases when $k\in[3,800]$. Note that for these values of the parameter $k$, Lemma \ref{cota_nl} gives us absolute upper bounds for $n$ and $m$. However, these upper bounds are so large that we wish to reduce them to a range where the solutions can be identified by using a computer. To do this, we let 
\begin{equation} \label{eq:z1}
z_1:=(n-1)\log\alpha-(m-1)\log\beta+\log \mu(k,\ell),
\end{equation}
where $\mu(k,\ell):=(2\alpha-1)f_k(\alpha)((2\beta-1)f_{\ell}(\beta))^{-1}$. Therefore, \eqref{deslambda} can be rewritten as 
\begin{equation}\label{dese^z1}
|e^{z_1}-1|<\frac{6}{\beta^{m-1}}.
\end{equation} 
Since $z_1\neq 0$ we distinguish the following cases. If $z_1>0$, then it follows from \eqref{dese^z1} that
\[
0<z_1\leq e^{z_1}-1<\frac{6}{\beta^{m-1}}.
\]
Replacing $z_1$ in the above inequality by its formula \eqref{eq:z1} and dividing both sides of the resulting inequality by $\log\beta$, we get	
\begin{equation}\label{z1>0-a}
0<(n-1)\left(\frac{\log\alpha}{\log\beta}\right)-m+\left(1+\frac{\log\mu(k,\ell)}{\log\beta}\right)<13\cdot\beta^{-(m-1)},
\end{equation}
where we have used the fact $1/\log\beta<2.1$ once again.  We put
\[
\hat{\gamma}:=\hat{\gamma}(k,\ell)=\frac{\log \alpha}{\log\beta},\quad \hat{\mu}:=\hat{\mu}(k,\ell)=1+\frac{\log\mu(k,\ell)}{\log\beta}, \quad A:=13,\quad {\text{\rm and}}\quad B:=B(\ell)=\beta.
\]
We also put $M_k:=\left\lfloor 5.4\times 10^{14} k^8 \log^3 k\right\rfloor$, which is an upper bound on $n$ by Lemma \ref{cota_nl}. The fact that $\hat{\gamma}$ is an irrational number can be found in \cite[p. 2129]{BL3}. Thus, the above inequality \eqref{z1>0-a} yields
\begin{equation} \label{z1>0-b}
0<(n-1)\hat{\gamma}-m+\hat{\mu}<A\cdot B^{-(m-1)}.
\end{equation}
It then follows from Lemma \ref{reduce}, applied to inequality \eqref{z1>0-b}, that 
\[
m-1<\frac{\log(Aq/\epsilon)}{\log B},
\] 
where $q=q(k,\ell)>6M_k$ is a denominator of a convergent of the continued fraction of $\hat{\gamma}$ such that $\epsilon=\epsilon(k,\ell)=||\hat{\mu}q||-M_k||\hat{\gamma}q||>0$. A computer search with \emph{Mathematica} revealed that if $k,\ell\in[2,800]$ with $\ell<k$, then the maximum value of $\log(Aq/\epsilon)/\log B$ is $<$ 1600. Hence, we deduce that the possible solutions $(n,k,m,\ell)$ of the equation \eqref{eqn} for which $k,\ell$ are in the range $[2,800]$ with $\ell<k$ and $z_1>0$, all have $m\in[7,1600]$. 

Next we treat the case $z_1<0$. First of all, one checks easily that $6/\beta^{m-1}<1/2$ for all $\ell\geq 2$ since $m\geq 7$. Thus, from \eqref{dese^z1}, we have that $|e^{z_1}-1|<1/2$ and therefore $e^{|z_1|}<2$. Since $z_1<0$, we have
\[
0<|z_1|\leq e^{|z_1|}-1=e^{|z_1|}|e^{z_1}-1|<\frac{12}{\beta^{m-1}}.
\]
In a similar way as in the case when $z_1>0$, and by recalling that $1/\log\alpha<2$ (since $k\geq 3$), we obtain
\begin{equation}\label{z_1<0}
0<(m-1)\hat{\gamma}-n+\hat{\mu}<A\cdot B^{-(m-1)},
\end{equation}
where now 
\[
\hat{\gamma}:=\hat{\gamma}(k,\ell)=\frac{\log \beta}{\log\alpha},\quad \hat{\mu}:=\hat{\mu}(k,\ell)=1-\frac{\log\mu(k,\ell)}{\log\alpha}, \quad A:=24,\quad {\text{\rm and}}\quad B:=B(\ell)=\beta.
\]
Here, we also took $M_k:=\left\lfloor 5.4\times 10^{14}\,k^8\,\log^3 k\right\rfloor$, which is an upper bound on $m$ by Lemma \ref{cota_nl}, and we applied Lemma \ref{reduce} to inequality \eqref{z_1<0} for each $k,\ell\in[2,800]$ with $\ell<k$.  In this case, with the help of \emph{Mathematica}, we found that the maximum value of $\log(Aq/\epsilon)/\log B$ is also $<$ 1600.  Thus, the possible solutions $(n,k,m,\ell)$ of the equation \eqref{eqn} for which $k,\ell$ are in the range $[2,800]$ with $\ell<k$ and $z_1<0$, all have $m\in[7,1600]$. 

Finally, we use \emph{Mathematica} to compare $L_n^{(k)}$ and $L_m^{(\ell)}$ for the range $6\leq n,m \leq 1600$ and $2\leq k,\ell\leq 800$, with $n<m$, $\ell<k$ and checked that the only solutions of the equation \eqref{eqn} in this range are the trivial solutions given by \eqref{trivial}.  This completes the analysis in the case $k\in[3,800]$.


\section{The case of large $k$}

\noindent From now on, we assume that $k>800$. For such $k$ we have
\[
n<m< 5.4\times 10^{14}\,k^8\,\log^3 k<2^{k/2}.
\]
It then follows from Lemma \ref{estimaciondeltaeta} that
\[
(2\alpha-1)f_k(\alpha)\alpha^{n-1}=3\cdot 2^{n-2}+3\cdot 2^{n-1}\eta_1 + \frac{\delta_1}{2}+\eta_1\,\delta_1,
\]
where $\eta_1$ and $\delta_1$ are real numbers such that
\[
|\eta_1|<\frac{2k}{2^{k}}  \quad \text{and} \quad  |\delta_1|<\frac{2^{n+2}}{2^{k/2}}.
\]
So, from the above equality, we get
\begin{equation}\label{(A1)}
\left|(2\alpha-1)f_k(\alpha)\alpha^{n-1}-3\cdot 2^{n-2}\right|<3\cdot \frac{2^{n}k}{2^{k}}+\frac{2^{n+1}}{2^{k/2}}+\frac{2^{n+3}k}{2^{3k/2}}<15\cdot \frac{2^{n-2}}{2^{k/2}},
\end{equation}
where the last inequality holds because $k>800$.  We will use estimate \eqref{(A1)} later. Let us now get some absolute upper bounds for the variables. In order to do so, we distinguish two cases.


\subsection{The case $m\leq 2^{\ell/2}$}

In this case, by using Lemma \ref{estimaciondeltaeta} once more, we get that
\[
(2\beta-1)f_{\ell}(\beta)\beta^{m-1}=3\cdot 2^{m-2}+3\cdot 2^{m-1}\eta_2 + \frac{\delta_2}{2}+\eta_2\,\delta_2,
\]
where now $\eta_2$ and $\delta_2$ are real numbers such that
\[
|\eta_2|<\frac{2\ell}{2^{\ell}}  \quad \text{and} \quad  |\delta_2|<\frac{2^{m+2}}{2^{\ell/2}}.
\]
Then the same argument used to derive \eqref{(A1)} leads to
\begin{equation}\label{(B1)}
\left|(2\beta-1)f_{\ell}(\beta)\beta^{m-1}-3\cdot 2^{m-2}\right|<45\cdot\frac{2^{m-2}}{2^{\ell/2}}
\end{equation}
for all $\ell \geq 2$. Hence, using \eqref{(A1)} and \eqref{(B1)}, we get
\[
\left|(3\cdot 2^{m-2}-3\cdot 2^{n-2})-((2\beta-1)f_{\ell}(\beta)\beta^{m-1}-(2\alpha
-1)f_{k}(\alpha)\alpha^{n-1})\right|<15\cdot \frac{2^{n-2}}{2^{k/2}}+45\cdot\frac{2^{m-2}}{2^{\ell/2}}
\]
giving
\[
\left|2^{m-2}-2^{n-2}\right|<19\cdot \frac{2^{m-2}}{2^{\ell/2}},
\]
where we used \eqref{des-al-beta} and the condition $\ell \leq m-1$. Dividing the last inequality above by $2^{m-2}$, we get 
\[
\frac{1}{2}\leq 1-2^{-(m-n)}<\frac{19}{2^{\ell/2}}.
\]
So, $2^{\ell/2}<38$ and therefore $\ell\leq 10$. Recalling that we are treating the case $m\leq 2^{\ell/2}$, it follows that $n<m\leq 37$. But a quick inspection of the list of generalized Lucas numbers tells us that the only solutions $(n,k,m,\ell)$ of equation \eqref{eqn} with $n\leq 36$, $k>800$, $m\leq 37$ and $\ell \leq 10$ are the trivial solutions given by \eqref{trivial}. This completes the analysis when $m\leq 2^{\ell/2}$.


\subsection{The case $2^{\ell/2}<m$}

Here, we have the following chain of inequalities
\[
2^{\ell/2}<m<5.4\times 10^{14}\,k^8\,\log^3 k<k^{14},
\]
which follow directly from Lemma \ref{cota_nl} together with the fact that $k>800$.  In particular,  
\begin{equation}\label{cotal(k)}
\ell<41\log k.
\end{equation}
On the other hand, combining \eqref{des-al-beta} and \eqref{(A1)}, we get 
\begin{align*}
\left|(2\beta-1)f_{\ell}(\beta)\beta^{m-1}-3\cdot 2^{n-2}\right|&<\left|(2\alpha-1)f_k(\alpha)\alpha^{n-1}-(2\beta-1)f_{\ell}(\beta)\beta^{m-1}\right|+15\cdot\frac{2^{n-2}}{2^{k/2}}\\
&<3+15\cdot\frac{2^{n-2}}{2^{k/2}}.
\end{align*}
Dividing both sides above by $3\cdot 2^{n-2}$, we arrive at
\[
\left|2^{-(n-2)}\cdot\beta^{m-1}\cdot 3^{-1}(2\beta-1)f_{\ell}(\beta)-1\right|<\frac{1}{2^{n-2}}+\frac{5}{2^{k/2}},
\]
which implies
\begin{equation}\label{lambda1}
\left|2^{-(n-2)}\cdot\beta^{m-1}\cdot 3^{-1}(2\beta-1)f_{\ell}(\beta)-1\right|<\frac{6}{2^{\Gamma}},
\end{equation}
where $\Gamma:=\min\{k/2,n-2\}$. The proof that the left--hand side of inequality \eqref{lambda1} is not zero is quite analogous to that given in \cite[p. 2134]{BL3}. We omit the details. 

We now lower bound the left--hand side of inequality \eqref{lambda1} using linear forms in logarithms. Here, Matveev's theorem together with a straightforward calculation, implies that
\begin{equation} \label{cotagamma}
\Gamma<4.46\times 10^{12}\,\ell^4\,\log^2\ell\,\log m. 
\end{equation}
Now, let $z_2:=(m-1)\log\beta-(n-2)\log 2+\log \mu(\ell)$ with $\mu(\ell)=3^{-1}(2\beta-1)f_{\ell}(\beta)$. So, estimation \eqref{lambda1} can be written as  
\begin{equation}\label{1-ez}
|e^{z_2}-1|<\frac{6}{2^{\Gamma}}.
\end{equation}
We distinguish two cases according to whether $z_2$ is positive or negative. First,  if $z_2>0$, then it follows from \eqref{1-ez} that
\[
0<z_2\leq e^{z_2}-1<\frac{6}{2^{\Gamma}}.
\]
Thus, 
\begin{equation}\label{z_2>0}
0<(m-1)\left(\frac{\log\beta}{\log 2}\right)-n+\left(2+\frac{\log{\mu(\ell)}}{\log 2}\right)<9\cdot 2^{-\Gamma}.
\end{equation}
We next treat the case $z_2<0$. First of all, observe that $6/2^{\Gamma}<1/2$ since $k>800$ and $n\geq 6$. Thus, $|e^{z_2}-1|<1/2$ leading to $e^{|z_2|}<2$. So, from \eqref{1-ez}, we get 
\[
0<|z_2|\leq e^{|z_2|}-1=e^{|z_2|}|e^{z_2}-1|<\frac{12}{2^{\Gamma}}.
\]
Consequently,
\begin{equation}\label{z_2<0}
0<(n-2)\left(\frac{\log 2}{\log\beta}\right)-m+\left(1-\frac{\log{\mu(\ell)}}{\log \beta}\right)<26\cdot 2^{-\Gamma}.
\end{equation}

In order to find some absolute upper bounds, we distinguish two subcases.


\subsubsection{Case 1. $\Gamma=k/2$}

Here, Lemma \ref{cota_nl}, together with bounds \eqref{cotal(k)} and \eqref{cotagamma}, yields
\[
k<2(4.46\times 10^{12})(41\log k)^4\,\log^2(41\log k)\log(5.4\times 10^{14}k^8\log^3 k).
\]
Using \emph{Mathematica} we obtained $k<2.8\times 10^{31}$. By Lemma \ref{cota_nl} once again and \eqref{cotal(k)}, we get $n<m<7.75\times 10^{271}$ and $\ell\leq 2970$. We record our conclusion as follows.

\begin{lemma}\label{cotas_absnlk}
If $(n,k,m,\ell)$ is a nontrivial solution in positive integers of equation \eqref{eqn} with $n\geq 6$, $k>800$, $2^{\ell/2}<m$ and $k/2\leq n-2$, then inequalities 
\[
n<m<7.75\times 10^{271},\quad k<2.8\times 10^{31}\quad\text{and}\quad\ell\leq 2970.
\]
hold.
\end{lemma}

We now reduce our previous bounds by using again Lemma \ref{reduce}. To avoid unnecessary repetitions, we consider only the case $z_2>0$. In this case, we take $M:=7.75\times10^{271}$ and we use Lemma \ref{reduce} on \eqref{z_2>0} for each $\ell\in[2,2970]$. A computer search with \emph{Mathematica} revealed that the maximum value of $k/2$ is at most 2980. Hence, we deduce that the possible solutions $(n,k,m,\ell)$ of the equation \eqref{eqn} for which $\ell\leq 2970$ and $z_2>0$ all have $k<5960$, and then from Lemma \ref{cota_nl} and \eqref{cotal(k)} we get $\ell\leq 360$ and $n<m<5.7\times 10^{47}$.

With this new upper bound for $m$ we repeated the process; i.e., we applied again Lemma \ref{reduce} with $M:=5.7\times 10^{47}$ for each $\ell\in[2,360]$. Here, we finally obtain that $k<740$, which is a contradiction. The same conclusion was obtained in the case $z_2<0$.


\subsubsection{Case 2. $\Gamma=n-2$}

We recall that we are in the situation $2^{\ell/2}<m$. Thus,   
\begin{equation}\label{cotal(m)}
\ell<\frac{2\log m}{\log 2}<3\log m.
\end{equation}
This, together with the bounds \eqref{util2} and \eqref{cotagamma}, tells us
\[
\frac{2m-11}{3}<4.46\times 10^{12}(3\log m)^4\,\log^2(3\log m)\,\log m.
\]
Using \emph{Mathematica}, we get an absolute upper bound for $m$, namely $m<9.1\times 10^{24}$. So, from \eqref{cotal(m)}, we get $\ell\leq 180$. We record what we have just proved.

\begin{lemma}\label{cotas_absnlk2}
If $(n,k,m,\ell)$ is a nontrivial solution in positive integers of equation \eqref{eqn} with $n\geq 6$, $k>800$, $2^{\ell/2}<m$ and $n-2<k/2$, then inequalities 
\[
n<m<9.1\times 10^{24}\quad\text{and}\quad\ell\leq 180
\]
hold.
\end{lemma}

Now, we would like to reduce our bound on $n$. If $z_2>0$, then we take $M:=9.1\times10^{24}$, which is an upper bound on $m$ from Lemma \ref{cotas_absnlk2}, and we use Lemma \ref{reduce} on inequality \eqref{z_2>0} for each $\ell\in[2,180]$.  

\emph{Mathematica} revealed that the maximum value of $n-2$ is at most 185. Hence, we deduce that the possible solutions $(n,k,m,\ell)$ of the equation \eqref{eqn} for which $\ell\leq 180$ and $z_2>0$ all have $n\leq 190$, and then from \eqref{util2} and \eqref{cotal(m)}, we get $m\leq 290$ and $\ell\leq 17$, respectively. The same conclusion remains valid for the case $z_2<0$.

Thus, we have reduced our problem to finding the solutions of \eqref{eqn} in the following range: $2\leq \ell \leq 17$, $6\leq n \leq 190$, $\ell+1<m \leq 290$ and $k>800$. But, for these values of $n$ and $k$, we have that $L_n^{(k)}=3\cdot 2^{n-2}$. Therefore, the problem is reduced to finding all the solutions of the equation
\begin{equation}\label{eqnfinal}
L_m^{(\ell)}=3\cdot 2^{n-2} \quad\text{with}\quad 2\leq \ell \leq 17,\quad \ell+1<m\leq 290 \quad \text{and}\quad  6\leq n\leq 190.
\end{equation}
Finally, a quick check with a computer confirms that equation \eqref{eqnfinal} has no solutions.  Thus, Theorem \ref{teo1} is proved.

\section{Acknowledgements}

\noindent We thank the  referee for a careful reading of the paper and for comments and suggestions which improved its quality. During the preparation of this paper, J. J. B. was partially supported by Universidad del Cauca and Colciencias from Colombia, and F. L. was supported in part by Project PAPIIT IN104512, UNAM, Mexico.

\end{document}